\documentclass[11pt,a4paper,reqno]{amsart} 
\usepackage{amsmath, amsthm, amssymb,  mathrsfs}
\usepackage{graphicx}
\usepackage{float}

\newcommand{\koniec}{\begin{flushright}  $\Box $ \end{flushright}}
\newtheorem{theo}{Theorem}[section]
\newtheorem{defi}[theo]{Definition}
\newtheorem{prop}[theo]{Proposition}  
\newtheorem{lemma}[theo]{Lemma}

\newcommand{\CP}{\mathbb{CP}}

\newcommand{\RP}{\mathbb{RP}}
\newcommand{\R}{\mathbb{R}}
\newcommand{\spp}{\mathbb{S}}
\def\z{{\mathcal T}}
\def\xx{{\mathcal X}}
\def\yy{{\mathcal Y}} 
\def\zz{{\mathcal Z}}
\def\p{\partial}

\def\e{{\bf e}}

\def\be{\begin{equation}}
\def\ee{\end{equation}}
\def\Th{\Theta}
\def\ep{{\varepsilon}}
\def\bea{\begin{eqnarray}}
\def\eea{\end{eqnarray}}

\def\l{\lambda}
\def\ov{\overline}

\newcounter{mnotecount}[section]

\renewcommand{\themnotecount}{\thesection.\arabic{mnotecount}}

\newcommand{\mnote}[1]
{\protect{\stepcounter{mnotecount}}$^{\mbox{\footnotesize
$
\bullet$\themnotecount}}$ \marginpar{
\raggedright\tiny\em
$\!\!\!\!\!\!\,\bullet$\themnotecount: #1} }

\date{September  27, 2012}

\begin{document}

\title{Twistor geometry of a pair of second order ODEs}

\author{Stephen Casey}
\address{Department of Applied Mathematics and Theoretical Physics\\ 
University of Cambridge\\ Wilberforce Road, Cambridge CB3 0WA, UK.
}
\email{sc581@cam.ac.uk}

\author{Maciej Dunajski}
\address{Department of Applied Mathematics and Theoretical Physics\\ 
University of Cambridge\\ Wilberforce Road, Cambridge CB3 0WA, UK.
}
\email{m.dunajski@damtp.cam.ac.uk}

\author{Paul Tod}
\address{The Mathematical Institute\\
Oxford University\\
24-29 St Giles, Oxford OX1 3LB\\ UK.
}
\email{tod@maths.ox.ac.uk}

\maketitle
\begin{center}
{\em Dedicated to Mike Eastwood on the occasion of his 60th birthday.}
\end{center}
\begin{abstract}
We discuss the twistor correspondence between path geometries in three dimensions with vanishing Wilczynski invariants and anti-self-dual conformal structures of signature $(2, 2)$. We show how to reconstruct a system of ODEs with vanishing invariants for a given conformal structure, highlighting the Ricci-flat case in particular. Using this framework, we give a new derivation of the Wilczynski invariants for a system of ODEs whose solution space is endowed with a conformal structure.  We explain how to reconstruct the conformal structure directly from 
the integral curves, and present new examples of systems of ODEs 
with point symmetry algebra of dimension four and greater which give rise to 
anti--self--dual structures with conformal symmetry algebra of the same dimension. Some of these examples are $(2, 2)$ analogues of plane wave space--times
in General Relativity. 
Finally we discuss a variational principle for twistor curves arising from
the Finsler structures with scalar flag curvature.
\end{abstract} 

\section{Introduction}
A path geometry on an open set $\z\subset \R^n$ 
consists of unparametrised curves, one
curve through each point in $\z$ in each direction. Regarding these curves as solutions
to a system of $(n-1)$ second order ODEs, one can give an alternative definition of the path
geometry as an equivalence class of systems of second order ODEs, where
two systems are regarded as equivalent if they can be mapped into each other
by a change of dependent and independent variables.
In this paper we shall study three--dimensional path geometries encoded into a system of two second order ODEs
\be
\label{system}
Y''=F(X, Y, Z, Y', Z'), \quad Z{''}=G(X, Y, Z, Y', Z'),
\ee
where $Y'=dY/dX$ etc, and $(F, G)$ are arbitrary functions on an open set in $\R^5$ which we assume to be of class $C^5$.
Two such systems are locally  equivalent if they are related by  
a point transformation
$
(X, Y, Z)\rightarrow (\ov{X}(X, Y, Z), \ov{Y}(X, Y, Z), \ov{Z}(X, Y, Z) )
$.
A natural question arises: Given a system (\ref{system}) is it equivalent
to a pair of trivial ODEs $Y''=0, \;\; Z''=0$?
The answer comes down to constructing a set of invariants
for (\ref{system}). These invariants can be divided into two groups ${\mathcal C}$ and ${\mathcal P}$. Vanishing of the invariants
from each of these groups selects a subclass of three--dimensional path geometries. 
\begin{itemize}
\item The system (\ref{system}) belongs to the {\em conformal branch} 
if its four--dimensional solution space $M$ admits a conformal structure
such that the two--dimensional surfaces in $M$ corresponding to points
in ${\mathcal T}$ are totally isotropic (this condition
uniquely determines the conformal structure).
The conformal branch is characterised by the vanishing of three fundamental
Wilczynski invariants ${\mathcal C}$ given by the expressions (\ref{wilczynski_sys}). 
Systems with vanishing Wilczynski invariants are also referred to as \emph{torsion-free}
\cite{Gro} and we will use this terminology several times throughout. 
\item The system (\ref{system}) belongs to the {\em projective branch} 
if its integral curves in $\mathcal{T}$ are unparametrised geodesics
of some projective connection.
The projective  branch is characterised by the vanishing of eight invariants
${\mathcal P}$ given by the expressions (\ref{fels}). 
\end{itemize}
A system is point equivalent to a trivial one if and only if it belongs to 
both the conformal and projective branches \cite{Gro}. In this case the symmetry algebra
of (\ref{system}) is isomorphic to $\mathfrak{sl}(4, \R)$. The corresponding Lie group $PSL(4, \R)$ acts projectively on $\z=\RP^3$ preserving the unparametrised geodesics of the flat projective connection. 

In this paper we shall concentrate on the conformal branch, and make use of the local isomorphism
\[
PSL(4, \R)\cong SO(3, 3),
\]
where $SO(3, 3)$ is the conformal group of the flat conformal structure on
the solution space $M$. If the curvature of the conformal structure does not vanish, then the conformal symmetry group is a proper subgroup of $SO(3, 3)$,
and conformal Killing vectors on $M$ give rise to point symmetries of the system (\ref{system}) (see Lemma \ref{lemma1}). We shall exploit this correspondence to construct
several examples of torsion-free ODE systems and the corresponding 
conformal structures which admit a large symmetry group, 9--dimensional
symmetry being the sub-maximal case.

In the rest of the paper we shall call ${\mathcal T}$
the {\it twistor space}, and the integral curves of (\ref{system})
the {\it twistor lines}. In the next section we shall give explicit expressions for
point invariants of (\ref{system}).
In section \ref{twist}
we summarise the twistor correspondence anti--self--dual conformal structures which underlies all constructions in this paper.
In section \ref{heavenly}, we describe how, given an ASD conformal structure of $(2, 2)$ signature, to construct a system of ODEs on the twistor space $\z$ whose integral curves are the twistor lines. Using this framework, in section \ref{paultod} we present a new derivation for the condition of vanishing Wilczynski invariants assuming solely that the space of solutions be endowed with a conformal structure. We demonstrate how this construction proceeds in the Ricci-flat case, and show (Theorem \ref{prop_paul}) that in this case the system of ODEs can be read-off directly from the solution of  Plebanski's
heavenly equation.  In section \ref{pointsymm}, we explore the local isomorphism between groups of point symmetries of (\ref{system})
and conformal symmetries of the corresponding  conformal structure.
In section \ref{sec_examples}
we construct examples of torsion-free ODE systems in three dimensions with point symmetry algebras of dimensions between
nine and two.  Finally in section \ref{section_finsler} we explore
the connection between systems of ODEs with vanishing Wilczynski invariants, and
unparametrised geodesics of Finsler structures with scalar flag curvature.
This correspondence gives rise to yet another class of examples.
Some of these examples have gravitational analogues in the theory of 
plane wave space--times, and the whole construction sheds new light
on variational aspects of the Nonlinear Graviton theorem.

\vskip5pt
It is a pleasure to
dedicate this paper to Mike Eastwood on the occasion of his 60th birthday.
Mike has made important contributions to twistor theory and conformal geometry over
the last thirty five years. He and his coworkers \cite{bceg} have developed
an approach to symmetry of differential equations which underlies some of our work.

\subsubsection*{Acknowledgements}
We are grateful to Boris Doubrov for telling us about the torsion--free
class of systems of ODEs, and for many useful discussions. In particular example
(\ref{boris}) was provided by Boris. We also thank Thomas Mettler for pointing out
the connection between the torsion--free systems and Finsler structures
with scalar flag curvature, and thank anonymous referees for pointing
out several inaccuracies in the original manuscript.

\section{Systems of ODEs and their invariants}
 \setcounter{equation}{0}
To present the invariants in a compact form set
$Y^A=(Y, Z)$ and $F^A=(F, G)$, where the capital letter indices 
$A, B, C, \dots $ take values $0, 1$. Moreover we set
\[
P^A={\frac{\p {Y^A}}{\p X}}, \quad{\mbox{and}}\quad
\frac{d}{d X}=\frac{\p}{\p X}+P^A\frac{\p}{\p Y^A}+F^A\frac{\p}{\p P^A}.
\]
The three Wilczynski conditions ${\mathcal C}$
are given by
\be
\label{wilczynski_sys}
T-\frac{1}{2}\mbox{Trace}(T){\mathbb{I}}=0, \quad\mbox{where}\quad
T^A_B=-\frac{\p F^A}{\p Y^B}-\frac{1}{4}\frac{\p F^A}{\p P^C}\frac{\p F^C}{\p P^B}
+\frac{1}{2}\frac{d}{d X}\frac{\p F^A}{\p P^B}.
\ee
Our terminology is motivated by the work of Doubrov \cite{Doubrov1,Doubrov}, who has demonstrated
that the conditions (\ref{wilczynski_sys}) for systems of non--linear
ODEs are equivalent to the classical Wilczynski invariants \cite{Wilczynski} on
the linearisations of these ODEs.

The complementary set of invariant conditions ${\mathcal P}$
characterising the projective branch is
\be
\label{fels}
{S^{A}}_{(BCD)}=0
\ee
where \cite{Fels}
\[
{S^A}_{BCD}=\frac{\p^3 F^A}{\p P^B\p P^C\p P^D}-\frac{3}{4}\frac{\p^3 F^E}{\p P^E\p P^B\p P^C}\delta^A_D.
\]
The following result gives an interpretation
of  the zero locus of Wilczynski invariants.
\begin{theo}[Grossman \cite{Gro}]
There is one--to--one correspondence between equivalence classes of
systems of second order ODEs with vanishing Wilczynski invariants (\ref{wilczynski_sys})
and conformal structures of signature $(2, 2)$ with vanishing self--dual
Weyl curvature.
\end{theo}
A point $p\in M$ corresponds to an integral curve $L_p$ of the system
(\ref{system}) -- a \emph{twistor line}. Moreover, if the Wilczynski invariants vanish then $M$ is endowed with a conformal structure and points in $\z$ correspond to special null two-surfaces in $M$ called \emph{$\alpha$--surfaces} (see next section). 
In Grossman's terminology $M$ acquires a \emph{Segre structure} from the torsion--free
system of ODEs \cite{Gro,Mettler}. In four dimensions the Segre structure is the splitting of the tangent
bundle into a product of rank two vector bundles. This is equivalent to the existence
of a conformal structure on $M$ - see formula (\ref{can_bun_iso}) in the next section -- so from now on we shall use the terminology of conformal geometry without referring to the Segre structures.

In the original Penrose twistor correspondence \cite{Pe76}, the four--manifold $M$ is assumed to be real analytic, and the correspondence extends to the complexified 
category. The curves $L_p$ in a complex three--fold $\z$ corresponding
to points in $M$ are globally characterised by the type of their normal bundle
$N(L_p)={\mathcal O}(1)\oplus{\mathcal O}(1)$, where ${\mathcal O}(k)\rightarrow \CP^1$
is a holomorphic line bundle with Chern class $k$. In the present context
we shall replace this global condition on $L_p$ by local--differential
conditions (\ref{wilczynski_sys}).

\section{Twistor Correspondence}
\label{twist}
Let $[g]=\{cg, c:M\rightarrow \R^+\}$ be a conformal structure on a 4-dimensional manifold $M$ consisting of an equivalence class of metrics of signature $(2, 2)$ and let $g\in [g]$. A vector field $V\subset TM$ is called null if $g(V, V)=0$. The null property of vector fields is invariant under rescalings of $g$, therefore it makes sense to talk about a null vector of a conformal structure. Conversely any conformal structure is completely characterised by specifying null vector fields.
\\Locally there exist real rank two vector bundles $\spp, \spp'$  (called spin-bundles) over $M$ equipped with parallel symplectic structures
$\ep, \ep'$ such that
\be
\label{can_bun_iso}
T M\cong {\spp}\otimes {\spp'}
\ee
is a  canonical bundle isomorphism, and
\[
g(v_1\otimes w_1,v_2\otimes w_2)
=\varepsilon(v_1,v_2)\varepsilon'(w_1, w_2)
\]
for $v_1, v_2\in \Gamma(\spp)$ and $w_1, w_2\in \Gamma(\spp')$.
Under this decomposition, any null vector field is of the form $V=\kappa\otimes\pi$ for some $\kappa\in\Gamma(\spp)$ and $\pi\in\Gamma(\spp')$. Also, the Riemann tensor can be decomposed as
\begin{eqnarray} \label{riemann}
R_{abcd} &=& \psi_{ABCD} \ep_{A'B'} \ep_{C'D'} + {\psi}_{A'B'C'D'}
\ep_{AB} \ep_{CD} \nonumber\\
&&+ \phi_{ABC'D'} \ep_{A'B'}\ep_{CD} +
\phi_{A'B'CD} \ep_{AB}\ep_{C'D'} \nonumber\\&& + \frac{R}{12} (\ep_{AC}\ep_{BD}
\ep_{A'C'} \ep_{B'D'} - \ep_{AD} \ep_{BC} \ep_{A'D'} \ep_{B'C'}),
\end{eqnarray}
where $\psi_{ABCD}$ and $\psi_{A'B'C'D'}$ are ASD and SD Weyl spinors
which are symmetric in their indices and
$\phi_{A'B'CD}=\phi_{(A'B')(CD)}$ is the traceless Ricci spinor.
\\An $\alpha$--plane is a two--dimensional plane in $T_pM$
spanned by null vectors of the above form with $\pi$ fixed, and
an $\alpha$--surface is a two--dimensional surface in $\zeta\subset M$ such that its tangent plane at every point is an $\alpha$--plane. Examining the Frobenius integrability conditions
for the existence of $\alpha$--surfaces leads to
the seminal result of Penrose \cite{Pe76}: A maximal, three dimensional, family of $\alpha$--surfaces exists
in $M$ iff the Weyl tensor of $g$ is anti--self--dual.
\\The anti--self--duality
of the Weyl tensor is a conformally invariant property, therefore one can talk about ASD conformal structures.
Let us assume that the conformal structure is real-analytic and therefore
can be complexified. Thus $(M, [g])$ is a holomorphic four--manifold
with a holomorphic conformal structure.
\begin{defi}
The twistor space ${\z}$
of a holomorphic conformal structure $(M, [g])$
with anti--self--dual Weyl curvature is the manifold
of $\alpha$--surfaces in $M$.
\end{defi}
The twistor correspondence puts  restrictions on the twistor space $\z$:
\begin{theo}[Penrose \cite{Pe76}]
\label{suptw1}
There is a one--to--one  correspondence between ASD conformal structures $(M, [g])$ and
three-dimensional complex
manifolds  ${\z}$ with the four parameter family of sections
of $\mu$ with normal bundle ${\mathcal O}(1)\oplus{\mathcal O}(1)$.
\end{theo}
The additional structures on the twistor space which give Ricci flat metric
$g\subset [g]$ are 
\begin{enumerate}
\item A projection $\mu :{\z}\longrightarrow \CP^1$,
such that the four parameter family of curves above
are sections of $\mu$.
\item Symplectic structure with values in ${\mathcal O}(2)$ on the fibers
of $\mu$.
\end{enumerate}
The points in the twistor space ${\z}$ correspond to totally null self--dual
two dimensional surfaces (which coincide with the $\alpha$-surfaces) in $M$, and
points in $M$ correspond to rational curves ( in the Ricci flat case these are 
sections of $\mu$) in ${\z}$.
The relation between ${\z}$ and $M$ is
best revealed by exploiting the double fibration
picture
\be
\label{doublefib11}
{M}\stackrel{r}\longleftarrow
{\mathcal F}\stackrel{s}\longrightarrow {\z}
\ee
where the correspondence space
\[
{\mathcal F}={\z}\times { M}|_{\zeta\in L_p}= {M}\times\CP^1
\]
can be identified with the projective primed spin bundle
$P(\spp')$. A point in ${\mathcal F}$ corresponds to a point in $p\in M$ together
with one $\alpha$--surface containing $p$.

The holomorphic curves $s(\CP^1_p)$ where
$\CP^1_p=r^{-1}p$, $p\in {M}$, have normal bundle $N={\mathcal
O}(1)\oplus{\mathcal O}(1)$. Two points $p_1$ and $p_2$ in $M$ are null separated iff the corresponding curves $L_{p_1}$ and $L_{p_2}$ intersect  at one point.

In the double fibration approach, the twistor space arises as a quotient
of ${\mathcal F}$ by a two--dimensional integrable distribution (called the Lax pair) defined
by two vector fields $L_0, L_1$ which at each point  of ${\mathcal F}$
are horizontal lifts to $P(\spp')$ of vectors spanning an $\alpha$ surface at a point in $M$. The existence of the Lax pair for an arbitrary conformal structure with
 vanishing self-dual Weyl curvature has the following consequence  \cite{mason, maciej}.
\begin{theo}
\label{tetradlaxprop}
Let $E_1, \dots, E_4$ be four real  vector fields on $M$
and let $e^1, \dots, e^4$ be the corresponding dual one--forms. The
conformal structure defined by
\[
g=e^1\odot {e}^2 - e^3\odot {e}^4
\]
is  ASD  if and only if there exists functions 
$f_0, f_1$ on $M \times \RP^1$ depending on $\lambda\in \RP^1$ 
such that the distribution 
\be 
\label{tetradlax0} L_0 = {E}_1 -  \lambda {E}_3 + f_0 \frac{\p}{\p \lambda}, 
\quad L_1 =
{E}_4 -  \lambda {E}_2 + f_1 \frac{\p}{\p \lambda}  \ee
is Frobenius integrable, that is, $[L_0, L_1] = 0$ modulo $L_0$ and $L_1$.
\end{theo}

The real $(2, 2)$ ASD conformal structures are obtained by introducing
an involution $\tau:\z\rightarrow \z$ on the twistor space given by complex conjugation.
The fixed points correspond to real $\alpha$--surfaces in $M$.
The involution acts on the  twistor lines,
thus giving rise to maps from $\CP^1$ to $\CP^1$
which swaps the lower and upper hemispheres
preserving the real equator. The moduli space of $\tau$--invariant twistor curves
is a real four--manifold.

\section{From ASD Conformal Structures to Systems of ODEs}
 \setcounter{equation}{0}
\label{heavenly}
In this Section we shall present an algorithm for reconstructing a system
(\ref{system}) of torsion--free ODEs from a given ASD conformal structure
$[g]$. Given $(M, [g])$ we shall construct a real projective line 
parametrised by $\lambda\in\RP^1$ 
of real $\alpha$--surfaces through a point in $p\in M$. This is a curve 
in $M\times \RP^1$ which also depends on the four coordinates of $p$ in $M$. 
We shall regard
this as a parametrised integral curve of a system of two second order ODEs. The system (\ref{system}) will arise by eliminating the coordinates of $p$ and removing the parametrisation.
To implement this procedure,  consider the Lax pair
(\ref{tetradlax0}) (see section \ref{twist})
of vector fields  $L_0$, $L_1$ on $M\times \RP^1$. Now proceed as follows
\begin{enumerate}
\item Find three independent functions $(\xx, \yy, \zz)$ on ${\mathcal F}=M\times\RP^1$
which satisfy
\[
L_0 f=0,\qquad L_1 f=0.
\]
These functions descend to the twistor space ${\z}$ where
they provide a local coordinate system $(X, Y, Z)$.
A point in $M$ corresponds to a curve in ${\z}$.
The pull back of the four--parameter family of  curves to ${\mathcal F}$ can be parametrised as
\[
\lambda\longrightarrow (\xx(\lambda, p), \yy(\lambda, p),
\zz(\lambda, p)),
\]
where $p=(w, z, x, y)$ are local coordinates on $M$.
\item Use the implicit function theorem to solve the equation $X=\xx(\lambda, p)$ for
$\lambda=\lambda(X, p)$. Then solve the relations
\[
Y=\yy,\quad Z=\zz, \quad Y'=\frac{\p \yy}{\p X},
\quad Z'=\frac{\p \zz}{\p X}
\]
to express $(w, z, x, y)$ as functions of $(Y, Z, Y', Z')$.
\item Differentiating once more, and substituting for
$(w, z, x, y)$ from above gives a pair of second order ODEs (\ref{system}).
\end{enumerate}
As an example of this construction, let us consider the ASD Ricci-flat metric which, in local coordinates may be written as
\begin{equation}
g = dw dx + dz dy - \Theta_{xx} dz^2 - \Theta_{yy} dw^2 + 2 \Theta_{xy} dw dz
\label{nkmetric}
\end{equation}
where $\Theta_x=\p_x\Theta$ etc, and 
$\Theta = \Theta(w,x,y,z)$ is a function satisfying the \emph{second heavenly equation}
\begin{equation}
\Theta_{xw} + \Theta_{yz} + \Theta_{xx} \Theta_{yy}- \Theta_{xy}^2 = 0.
\label{secondeq}
\end{equation}
It follows from the work of Pleba\'nski \cite{Pl75} that all $(2, 2)$ ASD Ricci flat
metrics locally arise from some solution to this equation.
The only non--vanishing part of curvature is given by the ASD Weyl spinor 
 which in the spin-frame adapted to the heavenly equation is 
a section of $\mbox{Sym}^4(\spp)\rightarrow M$ given by
\[
\psi_{ABCD}=\frac{\p^4\Theta}{\p x^A\p x^B\p x^C \p x^D}
\]
where the indices $A, B, \dots$ take values $0, 1$ and
$x^A=(y, -x)$. The Lax pair for this system may be written as
\begin{eqnarray}
\label{Lax2}
L_0&=&\p_y-\l(\p_w-\Th_{xy}\p_y+\Th_{yy} \p_x),\nonumber\\
L_1&=&\p_x+\l(\p_z+\Th_{xx}\p_y-\Th_{xy} \p_x).
\end{eqnarray}
A curve $L_p\subset {\z}$ corresponding to a point $p\in M$
is parametrised by choosing a two--dimensional fiber of
$\mu:{\z}\rightarrow\CP^1$ and defining
$(w, z)$ to be the coordinates of the initial point of the curve,
and $(y, -x)$ to be the tangent vector to the curve.
Thus the pulled back curve is
$
\lambda\rightarrow (X=\lambda, Y=\yy(\lambda, p), Z=\zz(\lambda, p))
$
where the functions $(\yy, \zz)$ admit the following expansion \cite{DM00}
\begin{eqnarray}
\label{expan}
\yy&=&w+\lambda y-\Theta_x\,\lambda^2+\Theta_z\,\l^3+\dots\\
\zz&=&z-\lambda x-\Theta_y\,\lambda^2 -\Theta_w\,\l^3+\dots\nonumber
\end{eqnarray}
and the higher order terms can be determined by recursion on successive powers of $\lambda$ from $L_A(\yy)=L_A(\zz)=0$, where $A=0, 1$.
\\ \\
{\bf Example 1.} Consider a solution to (\ref{secondeq})
\[
\Theta=\frac{1}{4}{y^4}.
\]
The corresponding system of ODEs is 
\begin{equation}
Y''=0 , \qquad Z''=-2(Y')^3.
\label{submax}
\end{equation}
The system has a nine--dimensional symmetry group (see formula (\ref{algebra1})), 
which is also the largest symmetry group of a non--trivial ASD conformal structure. The corresponding $pp$--wave metric
has constant Weyl curvature with only one non-vanishing component $\psi_{0000}=6$.
\\ \\
{\bf Example 2.}
This example, when analytically continued to the Riemannian
signature, is relevant in the theory of gravitational instantons.
The second heavenly equation
(\ref{secondeq}) with $\Th_z=0$ reduces to a wave equation
on a flat $(2+1)$--dimensional background.
Introduce $t:=\Th_x$ and perform a Legendre transform
\[
H(t, y, w):=tx(w, y, t)-\Th(w, y, x(w, y, t)).
\]
Then $x=H_t$, $\Theta_y=-H_y$ and (\ref{secondeq}) yields the wave
equation 
\be
\label{wave111}
H_{tw}+H_{yy}=0.
\ee
together with the metric
\begin{eqnarray}
\label{GHmetric111}
g&=&H_{tt}(\frac{1}{4} d y^2+ d w d t)-\frac{1}{H_{tt}}
(d z-\frac{H_{tt}}{2} d y+H_{ty} d w)^2\nonumber\\
&=&V(\frac{1}{4} d y^2+ d w d t)-V^{-1}( d z+A)^2,
\end{eqnarray}
where  $V=H_{tt}$ and $A=H_{ty} d w-(H_{tt}/2) d y$ satisfy the
monopole equation
$
*dV=dA
$
where $*$ is the Hodge operator on $\R^{2, 1}$ with its flat metric.
Thus (\ref{GHmetric111})
is an analytic continuation of the  Gibbons--Hawking metric \cite{GH}.
The vanishing of $\Theta_z$ implies that the whole series
(\ref{expan}) for
$\yy$ truncates at 2nd order. From (\ref{expan}), we get
\begin{eqnarray}
\label{GHcurves}
\yy&=&w+\lambda y-\lambda^2t,\\
\zz&=& z-\lambda H_t+\lambda^2H_y+\lambda^3H_w +...\;\nonumber
\end{eqnarray}
where $H=H(w, y, t)$, from which we can obtain the corresponding path geometry. An example with $H(w,y,t) = yt^2$ is
\begin{eqnarray}
\label{fourdexam}
Y''&=&X^{-1}Y'-\sqrt{(X^{-1}Y')^2-2X^{-1}Z'},\\
Z''&=&(X^{-1}Y'- \sqrt{(X^{-1}Y')^2-2X^{-1}Z'})^2/2.\nonumber
\end{eqnarray}
The potential in the Gibbons--Hawking metric is linear in the flat coordinates
on $\R^{2, 1}$.
\section{Wilczynski Invariants}
 \setcounter{equation}{0}
\label{paultod}
We shall now demonstrate a new way of deriving the Wilczynski 
invariants (\ref{wilczynski_sys}), from the double fibration picture (\ref{doublefib11}). This will provide a converse of the construction presented in the last section.
Our procedure is analogous to the recursive construction of the Wilczynski
invariants of a single $n$th order ODE \cite{DT}.

Each point in $\z$ corresponds to an $\alpha$-plane in $M$, and the functions 
$P^{A}=(dY/dX, dZ/dX)$ are null coordinates which are mutually orthogonal and thus $a_0 P^0 + a_1 P^1$ is null for arbitrary constants $a_0, a_1$. Furthermore, by differentiation, the one-form $a_0 dY^0+ a_1 dY^1$ is orthogonal 
to $a_0 dP^0+ a_1 dP^1$. 

In the derivation below we shall
regard $(Y^A, P^A, X)$ as coordinates on the five dimensional 
correspondence space ${\mathcal F}$ from the double fibration picture 
(\ref{doublefib11}),
and define
the degenerate metric $g$ on this five--dimensional space. 
We then demand that this
quadratic form Lie derives up to scale along the total derivative $d/d X$,
and so it gives a conformal structure on $M$.
The metric $g$  necessarily takes 
the form
\begin{equation}
g = \varepsilon_{AB} dY^A dP^B + \phi_{AB} dY^A dY^B
\label{confstruc}
\end{equation}
where $\phi_{AB} = \phi_{(AB)}$. The conformal structure of $M$ is 
invariant along the fibres of $s$ in the fibration (\ref{doublefib11}),
and therefore
\begin{equation*}
\frac{dg}{d X} = \Omega^2 g
\end{equation*}
for some function $\Omega$. Plugging in the expression (\ref{confstruc}) and comparing coefficients on both sides we obtain the equations
\begin{eqnarray}
2\phi_{AB} + \frac{\partial F^C}{\partial P^B} \varepsilon_{AC} &=& \Omega^2 \varepsilon_{AB} \label{confeqone} \\
\frac{d\phi_{AB}}{d X} + \frac{\partial F^C}{\partial Y^{(B}} \varepsilon_{A)C} &=& \Omega^2 \phi_{AB} \label{confeqtwo}.
\end{eqnarray}
Taking the trace of the (\ref{confeqone}) (using $\varepsilon$ to raise and lower indices)
\begin{equation*}
\frac{\partial F^C}{\partial P^C} = 2 \Omega^2
\end{equation*}
and substituting back into (\ref{confeqone}) for $\Omega$ we get
\begin{equation*}
\phi_{AB} = -\frac{1}{2} \varepsilon_{AC} \frac{\partial F^C}{\partial P^B} + \frac{1}{4} \varepsilon_{AB} \frac{\partial F^C}{\partial P^C}.
\end{equation*}
Then, from (\ref{confeqtwo}), we obtain the Wilczynski invariants
(\ref{wilczynski_sys})
\begin{equation*}
-\frac{1}{2}\frac{d}{dX}\frac{\partial F^A}{\partial P^B} + \frac{\partial F^A}{\partial Y^B} + \frac{1}{4} \frac{\partial F^C}{\partial P^B} \frac{\partial F^A}{\partial P^C} \sim \delta^A_B
\end{equation*}
where here we have used the fact that
\begin{equation*}
\frac{\partial F^C}{\partial P^B} \frac{\partial F^A}{\partial P^C} - \frac{\partial F^C}{\partial P^C} \frac{\partial F^A}{\partial P^B} \sim \delta^A_B
\end{equation*}
to rewrite the third term.
\\Note here that the expression for the conformal metric (\ref{confstruc}) actually gives a metric if $\Omega = 0$ i.e, $\frac{\partial F^C}{\partial P^C} = 0$ which implies
\begin{equation}\label{pt1}
F^A = 2 \varepsilon^{AB} \frac{\partial \Lambda}{\partial P^B}
\end{equation}
for some function $\Lambda$. The metric (\ref{confstruc}) then resembles the heavenly 
form (\ref{nkmetric}) of the Ricci--flat metric. The exact equivalence
arises from evaluating (\ref{confstruc}) at $X=0$, where $Y^A=(w, z), P^A=(-y, x)$ 
and $\Lambda(Y^A, P^A)=-\Theta(w, z, x, y)$.
\\ \\
{\bf Example 3.} One example is given by
\begin{equation}
Y'' = 0 \quad Z''=\beta (Y')
\label{ptexample}
\end{equation}
for some arbitrary function $\beta$. The ASD conformal structure on the solution space $M$ is type--N and Ricci-flat. Notice here that the system (\ref{submax}) of sub-maximal point symmetry is just a special case of this solution.
\vskip5pt

The last example was Ricci-flat. One may seek conditions on $F^A$ for this to be true more generally. In other words one seeks to obtain (\ref{secondeq}) as well as (\ref{nkmetric}). For this, 
one may proceed as follows: introduce the two-form
\[\Sigma=dY\wedge dZ=\frac12\varepsilon_{AB}dY^A\wedge dY^B,\]
and then calculate
\[\frac{d\Sigma}{dX}=\varepsilon_{AB}dY^A\wedge dP^B,\]
\[\frac{d^2\Sigma}{dX^2}=\varepsilon_{AB}(dP^A\wedge dP^B+dY^A\wedge dF^B),\]
and
\[\frac{d^3\Sigma}{dX^3}=\varepsilon_{AB}(3dP^A\wedge dF^B+dY^A\wedge \frac{d}{dX}dF^B),\]
all modulo $dX$ and where $d/dX$ may be thought of as Lie-derivative. 

In the context of the Nonlinear Graviton Theorem \ref{suptw1} the two-form
$\Sigma$ is the pull back of the 
${\mathcal{O}}(2)$--valued symplectic structure on the fibres
of ${\z}\longrightarrow \CP^1$ to the correspondence space. Thus, up to 
an overall scale, this two--form
is given by a quadratic polynomial in the base coordinate $X$.
We can therefore examine the consequences of imposing the requirement
\[\frac{d^3\Sigma}{dX^3}\sim\Sigma\mbox{   modulo   }dX.\]
The coefficient of $dP^0\wedge dP^1$ on the left-hand-side must vanish, which imples $\partial F^B/\partial P^B=0$ and we may introduce $\Lambda$ as in (\ref{pt1}). Then the vanishing of terms in $dY^A\wedge dP^B$ reduces to a single equation 
\[\varepsilon^{AB}\frac{\partial^2\Lambda}{\partial Y^A\partial P^B}+\frac12\varepsilon^{AD}\varepsilon^{BE}\frac{\partial^2\Lambda}{\partial P^A\partial P^E}\frac{\partial^2\Lambda}{\partial P^B\partial P^D}=0,\]
which is precisely (\ref{secondeq}) in this notation. 
The remaining condition (\ref{confeqtwo}) can now be viewed as the 
initial value problem for $\Lambda$ with the initial data given by the solution of the heavenly equation. This proves the following
\begin{theo}
\label{prop_paul}
Let $\Theta=\Theta(w, z, x, y)$ be a solution to the heavenly equation
{\em(\ref{secondeq})} which gives the Ricci--flat ASD metric {\em(\ref{nkmetric})}.
The corresponding system of ODEs with vanishing Wilczynski invariants
is
\be
Y''=2\frac{\p\Lambda}{\p Z'}, \quad Z''=-2\frac{\p \Lambda}{\p Y'},
\ee
where $\Lambda|_{X=X_0}=-\Theta(Y(X_0), Z(X_0), Z'(X_0), -Y'(X_0))$
and the $X$--dependence of $\Lambda$ is determined by {\em(\ref{confeqtwo})}.
\end{theo}
The formula (\ref{confstruc}) gives a way to reconstruct the conformal 
structure directly from the system of ODEs. An alternative procedure
based on the original non--linear graviton construction of Penrose
\cite{Pe76} gives the conformal structure directly from the integral curves
of the system of ODEs. This procedure will be implemented in the example
(\ref{boris}).
\section{Symmetries of Torsion-Free Path Geometries}
 \setcounter{equation}{0}
\label{pointsymm}
In the twistor approach, we observed there is a one-to-one correspondence between points in $\z$ and $\alpha$-surfaces in $M$. Therefore, any transformation
which preserves points in $\z$ will give rise to a transformation of $M$ which preserves the ASD conformal structure. 
\begin{lemma}
\label{lemma1}
There is a one--to--one correspondence between conformal Killing vectors 
of a $(2, 2)$ ASD conformal structure  $(M, [g])$ and point symmetries of 
the torsion--free system of ODEs whose integral curves
are twistor lines for  $(M, [g])$.
\end{lemma}
{\bf Proof.} 
We shall use the double fibration picture (\ref{doublefib11}). 
Given a conformal 
Killing vector $K$ of $(M, [g])$ we can lift it to the vector field
$\widetilde{K}$  on the correspondence space
${\mathcal F}$, so that $[L_0, \widetilde{K}]=0$ and $[L_1, \widetilde{K}]=0$, where
$L_0, L_1$ is the twistor distribution (\ref{tetradlax0}), and
the commutators vanish modulo a linear combination of $L_0$ and $L_1$.
The lift is explicitly given by $\widetilde{K}=K+Q\p_\lambda$, where $Q$ is a
quadratic polynomial in $\lambda$ with coefficients depending on coordinates
on $M$. In the spinor notation $Q=\nabla_{AA'}{K^{A}}_{B'}\pi^{A'}\pi^{B'}$,
where $\pi^{A'}=(1, -\lambda)$.

The space $\z$ is a quotient if ${\mathcal F}$ by the distribution
spanned by $L_0$ and $L_1$ and thus ${\mathcal K} =s_* \widetilde{K}$
is a vector field on $\z$. Therefore it generates a one--parameter
group of transformations of $\z$ which therefore 
takes $\alpha$--surfaces in $M$ to $\alpha$--surfaces. As $K$ generates
diffeomorphisms of $M$ and integral curves of the system (\ref{system})
in $\z$ correspond to points in $M$, then the action generated by 
${\mathcal K}$ preserves the integral curves of (\ref{system}).
Thus it is a point symmetry on (\ref{system}).

Conversely, a point symmetry of (\ref{system}) corresponds to 
a transformation of $M$ which maps $\alpha$--surfaces to $\alpha$--surfaces.
Therefore it gives a conformal Killing vector on $M$.

\koniec
In the trivial case ($Y'' = 0, Z'' = 0$), the path geometry has point symmetry group $SL(4, \R)$ which is isomorphic to the symmetry group of the flat conformal structure $SO(3, 3)$. The sub-maximal case of a torsion-free system of ODEs with nine-dimensional point symmetry algebra (\ref{submax}) which corresponds to a Ricci-flat ASD conformal structure with only one constant non-vanishing component of the ASD Weyl tensor. 
Therefore the `gap' in $(2, 2)$ conformal geometry equals to $6=15-9$. This, 
by Lemma \ref{lemma1} coincides with the gap of path geometries 
in  three dimensions,
as it is known that the sub-maximal symmetry in the projective class ${\mathcal P}$
is only eight dimensional. This follows from the general result of  
Egorov \cite{eg} that all sub-maximally symmetric projective connections on $\R^n$ are transitive and have symmetry algebra of dimension $n^2-2n+5$. 
See \cite{krug} for an account
of the general theory of such `gaps'.

Example (\ref{submax}) is the unique (up to diffeomorphism) torsion-free path geometry with point symmetry 
algebra of sub-maximal size and we expect systems with $6, 7, 8$ symmetries to be also comparatively rare. To see explicitly why this occurs, consider the lift of the integral curves of an arbitrary path geometry (\ref{system}) to the second jet bundle $J^2(\z, \mathbb{R})$ which is a seven-dimensional manifold with local coordinates given by 
$\left(X, Y, Z, Y', Z', Y'', Z''\right)$. 
Any point symmetry of (\ref{system}) can locally be described by some vector field $\chi$ on $\z$ which we can prolong to a vector field $\mbox{pr}^{(2)}{\chi}$ over some open set in $J^2(\z,\mathbb{R})$. Then, the functions
\begin{equation*}
\Delta^1 = Y''-F(X, Y, Z, Y', Z'), \quad \Delta^2=Z''-G(X, Y, Z, Y', Z')
\end{equation*}
are constant along $\mbox{pr}^{(2)}{\chi}$.
\\Now suppose that the given path geometry admits a Lie point symmetry algebra of dimension five, generated by five vector fields which prolong to an integrable distribution $\tilde{L}$  of $J^2(\z,\mathbb{R})$ on rank at most five. In the generic case, when the rank is five, 
these vector fields will lie in the tangent bundle of some five-dimensional submanifold $U \subset J^2(\z,\mathbb{R})$. Given that the codimension of $U$ is 2, we can  construct two functions which $\Delta^0, \Delta^1$ which are invariant with respect to the action of $\tilde{L}$, i.e, given a five-dimensional Lie algebra $L$ of vector fields over $\z$, we can find a (non-trivial) path geometry (\ref{system}) which has $L$ as its point symmetry algebra. This statement is true subject to some 
regularity conditions. For example it is sufficient to demand that the  
first prolongation of the vector fields
to the first jet bundle $J^1(\z, \R)$ must not be contained
within some four-dimensional submanifold of the tangent bundle 
$T(J^1(\z, \R))$.

\vskip 5pt
Lie algebras $L$ of dimension six or greater will generically not give
rise to functions  $\Delta^0, \Delta^1$ , 
and there will be a constraint on finding non-trivial path geometries with point symmetry algebras of this size. In particular, we must require that the prolonged algebra forms a distribution of rank lower than six. If additionally, we impose the constraint that the Wilczynski conditions (\ref{wilczynski_sys}) hold, then this sudden decline of examples will be observed sooner (at dimension four rather than six). In the next section, we outline the prolongation procedure and present some examples of path geometries with 
$4, 5, 6, 7, 8, 9$ point symmetries together with some details of the Lie algebraic structure. Where possible, we also make use of the twistor approach to say something about the corresponding conformal structure on the space of solutions.
\section{Examples}
\label{sec_examples}
For a given path geometry in three dimensions (\ref{system}) the point symmetries may be found by the following well known procedure:
\begin{enumerate}
\item Let the generator of a point symmetry be
\begin{equation*}
\chi = \chi_{-1} \partial_{X} + \chi_0 \partial_{Y} + \chi_1 \partial_{Z}
\end{equation*}
where the $\chi_i$ are functions of $(X, Y, Z)$ for which we must solve.
\item Determine the first and second prolongations by
\begin{equation*}
\eta^{(1)}_A = \frac{d\chi_A}{d X} - P^A \frac{d\chi_{-1}}{d X}, \quad  \eta^{(2)}_A = \frac{d \eta^{(1)}_A}{d X} - Q^A \frac{d\chi_{-1}}{d X}
\end{equation*}
$A=0,1$, where $Q^A := \frac{d^2 Y^A}{d X^2}$.
\item The prolongation of the vector field $\chi$ to the jet bundle $J^2(\z,\mathbb{R})$ is given by
    \begin{equation*}
     \text{pr}^{(2)}(\chi) = \chi + \displaystyle \sum_{B=0}^1 \eta^{(1)}_B \frac{\partial}{\partial P^B} + \sum_{B=0}^1 \eta^{(2)}_B \frac{\partial}{\partial Q^B}.
     \end{equation*}
\item Then we determine point symmetries by finding functions $\chi_i$ from
\begin{equation*}
\text{pr}^{(2)}(\chi)(\Delta^0) \vert_{\Delta^0 = \Delta^1 = 0} = \text{pr}^{(2)}(\chi)(\Delta^1) \vert_{\Delta^0 = \Delta^1 = 0} = 0.
\end{equation*}
\end{enumerate}
\subsection{Special Ricci-Flat Case}
Consider the system (\ref{ptexample}) with $\beta$ real analytic. The sub-maximal torsion-free system (\ref{submax}) lies in this class so we might expect it to yield more examples of systems with a high number of point symmetries. The corresponding ASD conformal structure on the moduli space of solutions is given by
\begin{equation*}
g = dw dx + dz dy - \frac{1}{2} \beta'(y) dw^2
\end{equation*}
and is Ricci-flat, as it corresponds to (\ref{nkmetric}) 
with $\Theta=(1/2)\int \beta(y)dy$.
\\Under these conditions, it is not difficult to show that the point symmetry algebra contains a six-dimensional subalgebra $L_6 \subset \mathfrak{sl}(4, \R)$. It transpires that $L_6$ is solvable and in terms of point symmetries we may write it down explicitly
\begin{equation*}
L_6 = \text{span} \{\e_1 = \partial_X, \e_2 = \partial_Y, \e_3 = \partial_Z, \e_4 = X\partial_{Y}, \e_5 = Z \partial_{Y}, \e_6 = X \partial_X 
+2Y \partial_{Y} +Z \partial_{Z} \}.
\end{equation*}
However, as we have seen, there will be some special cases of (\ref{ptexample}) for which the point symmetry algebra is larger but will contain $L_6$ as a subalgebra.

\begin{prop}
Consider a system of two second order ODEs of the form (\ref{ptexample})
for some function $\beta$ of the form
\begin{equation*}
\beta(Y') = \displaystyle \sum_{k=0}^{\infty} \xi_k {(Y')}^k.
\end{equation*}
If $\beta$ is a quadratic function then the system (\ref{ptexample}) is diffeomorphic to a trivial one (and the symmetry group is 15-dimensional).
Otherwise, the symmetry algebra has dimension
{\em \begin{equation*}
12 - \text{Rank}\; (M_1) - \text{Rank}\; (M_2)
\end{equation*}\em}
where $M_1$ is a matrix with rows
\begin{equation*}
( \xi_k \,\,\,\, \xi_{k+1}) \,\,\,\,\,,\,\,\,\,\, k \geq 3
\end{equation*}
and $M_2$ is a matrix with rows
\begin{equation*}
( \xi_k \,\,\,\, (k-2) \xi_k \,\,\,\, (k-3) \xi_{k-1} \,\,\,\, -(k+1) \xi_{k+1}) \,\,\,\,\,,\,\,\,\,\, k \geq 3.
\end{equation*}
\end{prop}
{\bf Proof.} 
Without loss of generality, let us simplify the problem by making the diffeomorphism
\begin{equation*}
Z \rightarrow Z + \frac{1}{2} \xi_0 X^2 + \xi_1 X Y + \frac{1}{2} \xi_2 Y^2
\end{equation*}
so that we obtain the system
\begin{equation*}
Y'' = 0 \,\,\,,\,\,\, Z'' = \displaystyle \sum_{k=3}^{\infty} \xi_k (Y')^k
\end{equation*}
which has the same number of point symmetries as the original.
\\For a given vector field,
\begin{equation*}
\chi = \chi_{-1} \frac{\partial}{\partial X} + \chi_0 \frac{\partial}{\partial Y} + \chi_1 \frac{\partial}{\partial Z},
\end{equation*}
the expressions
\begin{equation*}
\text{pr}^{(2)} (\chi) (\Delta^0)|_{\Delta^0=\Delta^1=0} \,\,\, \text{and} \,\,\,  \text{pr}^{(2)} (\chi) (\Delta^1)|_{\Delta^0=\Delta^1=0}
\end{equation*}
are real analytic in $P^0$ and $P^1$ with coefficients which are functions of $X$, $Y$ and $Z$. For $\chi$ to be a symmetry of the system (4.5), each of these coefficients must vanish separately. This leads to the following system of differential equations in $\chi_{-1}$, $\chi_0$ and $\chi_1$:
\begin{equation*}
\frac{\partial^2 \chi_0}{\partial X^2} = \frac{\partial^2 \chi_1}{\partial X^2} = \frac{\partial^2 \chi_0}{\partial X \partial Z} = \frac{\partial^2 \chi_1}{\partial X \partial Y} = \frac{\partial^2 \chi_{-1}}{\partial Y \partial Z} = \frac{\partial^2 \chi_{-1}}{\partial Z^2} = \frac{\partial^2 \chi_0}{\partial Z^2} = 0,
\end{equation*}
\begin{equation*}
2 \frac{\partial^2 \chi_0}{\partial X \partial Y} = \frac{\partial^2 \chi_{-1}}{\partial X^2} \,\,\,,\,\,\, \frac{\partial^2 \chi_0}{\partial Y^2} = 2 \frac{\partial^2 \chi_{-1}}{\partial X \partial Y} \,\,\,,\,\,\, \frac{\partial^2 \chi_0}{\partial Y \partial Z} = \frac{\partial^2 \chi_{-1}}{\partial X \partial Z},
\end{equation*}
\begin{equation*}
 2 \frac{\partial^2 \chi_1}{\partial X \partial Z} = \frac{\partial^2 \chi_{-1}}{\partial X^2} \,\,\,,\,\,\, \frac{\partial^2 \chi_1}{\partial Y \partial Z} = \frac{\partial^2 \chi_{-1}}{\partial X \partial Y} \,\,\,,\,\,\, \frac{\partial^2 \chi_1}{\partial Z^2} = 2 \frac{\partial^2 \chi_{-1}}{\partial X \partial Z},
\end{equation*}
\begin{equation*}
\frac{\partial^2 \chi_1}{\partial Y^2} = 3 \xi_3 \frac{\partial \chi_0}{\partial X}\,\,\,,\,\,\, \xi_3 \frac{\partial \chi_0}{\partial Z} = \frac{\partial^2 \chi_{-1}}{\partial Y^2} \,\,\,,\,\,\,  \frac{\partial^2 \chi_{-1}}{\partial Y^2} + 3 \xi_3  \frac{\partial \chi_0}{\partial Z} = 0,
\end{equation*}
and, for $k \geq 3$,
\begin{equation*}
(k-3) \xi_k  \frac{\partial \chi_{-1}}{\partial Z} - (k+1) \xi_{k+1} \frac{\partial \chi_0}{\partial Z} \,\,\,,\,\,\, \xi_{k+1}  \frac{\partial \chi_0}{\partial Z} = \xi_k  \frac{\partial \chi_{-1}}{\partial Z},
\end{equation*}
\begin{equation*}
(k-2)\xi_k  \frac{\partial \chi_{-1}}{\partial X}+ \xi_k  \frac{\partial \chi_1}{\partial Z} + (k-3) \xi_{k-1}  \frac{\partial \chi_{-1}}{\partial Y} - (k+1) \xi_{k+1}  \frac{\partial \chi_0}{\partial X} -k \xi_k  \frac{\partial \chi_0}{\partial Y} = 0.
\end{equation*}
Enforcing the system to be non-trivial (i.e, not all $\xi_k = 0$) we obtain
\begin{equation*}
 \frac{\partial \chi_{-1}}{\partial Z} =  \frac{\partial \chi_0}{\partial Z} = \frac{\partial^2 \chi_{-1}}{\partial Y^2} = 0.
\end{equation*}
Then, the most general solution for $\chi$ is
\begin{eqnarray*}
\chi_{-1} &=& a_1 X^2 + a_2 XY + (a_3 + a_4) X + a_5 Y + a_6 \\
\chi_0 &=& a_1 XY + a_7 X + a_2 Y^2 + 2 a_3 Y + a_8 \\
\chi_1 &=& a_1 XZ + a_9 X + a_2 YZ + (a_3 + a_{10}) Z + \frac{\xi_3}{2} a_1 Y^2 + \left( \frac{3 \xi_3}{2} a_7 + a_{11} \right) Y + a_{12}
\end{eqnarray*}
where, for all $k \geq 3$,
\begin{equation*}
\xi_k a_2 + \xi_{k+1} a_1 = 0
\end{equation*}
and
\begin{equation*}
\xi_k a_{10} + (k-2) \xi_k a_4 +(k-3) \xi_{k-1} a_5 - (k+1) \xi_{k+1} a_7 - (k+1) \xi_k a_3 = 0
\end{equation*}
and the result follows.
\vskip2pt

Thus, in the non-trivial case, the point symmetry algebra of 
(\ref{ptexample}) has at most dimension 9. A little more work shows us that 
systems of dimension 8 are not attainable for this example (if the rank of $M_2$ is 1, then the rank of $M_2$ is 2 and if the rank of $M_1$ is 2 then the rank of $M_2$ is at least 3).
\koniec
\subsubsection*{Dimension 9}
A system with point symmetry algebra of dimension 9 is given by
(\ref{submax})
\begin{equation*}
Y'' = 0, \quad Z'' = -2(Y')^3.
\end{equation*}
A quick check of the expressions for $M_1$ and $M_2$ will show that the system
\begin{equation*}
Y'' = 0, \quad Z'' = \frac{(Y')^3}{1-Y'}
\end{equation*}
also admits a nine-dimensional point symmetry algebra but these systems can be shown to be diffeomorphic. The associated Lie algebra in the first case can be written as
\begin{equation}
\label{algebra1}
L_9 = L_6 \oplus \text{span} \{\e_7, \e_8, \e_9 \}
\end{equation}
with
\begin{equation*}
\e_7 = 3Y \partial_{Y} + Z \partial_{Z}\,\,,\,\, \e_8 = \frac{3}{2} Z^2 \partial_{Y} + X \partial_{Z}\,\,,\,\, \e_9 = \frac{1}{2} X^2 \partial_{X} + \left( \frac{1}{2} X Y + \frac{1}{4} Z^3 \right) \partial_{Y} + \frac{1}{2} X Z \partial_{Z}.
\end{equation*}
This has Levi-decomposition
\begin{equation*}
L_9 = \tilde{L}_6 \ltimes \mathfrak{sl}(2, \R) = \text{span} \{ \e_2, \e_3, \e_4, \e_5, \e_7, \e_8 \} \ltimes \text{span} \{ \e_1, \e_6 - \frac{1}{2} \e_7, \e_9 \}.
\end{equation*}
\subsubsection*{Dimension 7}
A path geometry with point symmetry algebra of dimension 7 is given by
\begin{equation*}
Y'' = 0, \quad Z''= (Y')^k
\end{equation*}
for $k \geq 4$. This will be a solvable Lie algebra for any value of $k$ and is obtained by adding the vector field
$
\e_7 = k Y \partial_{Y} + Z \partial_{Z}
$
to $L_6$. 
\subsection{ASD Einstein - dimension 8}
Consider the system
\be
\label{boris}
Y'' = 0, \quad Z'' = \frac{2(Z')^2 Y'}{Z Y'-1}.
\ee
This system has been suggested to us by Boris Doubrov, who has also 
pointed out its  appearance in the theory of chains in the
homogeneous contact geometry \cite{cap}. The Wilczynski invariants
vanish for (\ref{boris}).

The corresponding conformal structure is
found by demanding that two neighbouring integral curves
in $\z$ intersect at one point  -- this is the condition
which selects null vectors in $M$. This, in the complexified setting, 
reduces the normal bundle
of the integral curve $L_p$ to ${\mathcal O}(1)\oplus{\mathcal O}(1)$.
To perform this calculation explicitly, solve (\ref{boris})
for $(Y, Z)$ and demand that equations $\delta Y=0$ and  $\delta Z=0$
for $X$ have a common solution. We apply this procedure to the integral curves
\[
\yy=w+yX,\quad \zz=\frac{1}{y}+\frac{1}{y^2(z-xX)}
\]
where $\delta Y=\sum_a\p_a \yy\; \delta x^a , 
\delta Z=\sum_a\p_a \zz \;\delta x^a$,
and $x^a=(w, z, x, y)$. 
Then the discriminant condition for two curves to intersect gives the metric
\begin{equation*}
g = dw dx + dz dy + x^2 dw^2 + \left( z^2 + \frac{2z}{y} \right) dy^2 + 2\left( zx + \frac{x}{y} \right) dw dy
\end{equation*}
which is ASD and Einstein, with scalar curvature equal to $-24$. Both the system of ODEs and the metric have eight-dimensional symmetry group $SL(3, \mathbb{R})$ which acts isometrically on $M$. The homogeneous model for the space of solutions is $SL(3, \mathbb{R})/GL(2, \mathbb{R})$, which is a $(2,2)$ real form of the Fubini-Study metric on $\mathbb{C}\mathbb{P}^2 = SU(3)/U(2)$.

\subsection{Symmetry algebra of dimension 5}
The example discussed in the Gibbons-Hawking context (\ref{fourdexam})
has five-dimensional point symmetry algebra which we can write as 
\begin{equation*}
{L_5} = \text{span} \{ \partial_{Y}, \partial_{Z}, X \partial_{X} + 
Y \partial_{Y}, -\frac{1}{2}X \partial_{X} + Z \partial_{Z}, 
X^2\p_Y+2Y\p_Z\}.
\end{equation*}
The algebra $L_5$ is solvable and  contains 
both Bianchi II and Bianchi V as three--dimensional subalgebras.
\subsubsection{Sparling-Tod solution}
Another way to construct a torsion-free path geometry with a given point symmetry algebra is to determine an ASD metric 
with that conformal symmetry algebra and determine the system of ODEs governing twistor lines via the method given in Section \ref{paultod}.
The metric of \cite{sparlingtod} given by
\begin{equation*}
g= dw dx + dz dy - \frac{2}{(xw+yz)^3} \left(w dz - z dw \right)^2
\end{equation*}
is an ASD Ricci-flat metric (\ref{nkmetric}) with $\Theta = \frac{1}{xw+yz}$ and has five-dimensional conformal symmetry algebra. 
The conformal  Killing vectors are
\[
{\bf k}_1= z\p_w - x\p_y,
{\bf k}_2= -2w\p_w + x\p_x - y\p_y ,
{\bf k}_3= w\p_w + z\p_z,
{\bf k}_4= -y\p_x + w\p_z,
{\bf k}_5= z\p_x - w\p_y.
\]
We should note here that the resulting system of ODEs does not coincide with the previous example as the conformal symmetry algebra of the Sparling-Tod solution is not solvable.
It has an $\mathfrak{sl}(2, \R)$ subalgebra generated by 
$({\bf k}_1, {\bf k}_2+ {\bf k}_3,
{\bf k}_4)$. Its Levi decomposition can be expressed as the semi-direct product of an $\mathfrak{sl}(2, \R)$  with the 2--dimensional
 non-Abelian Lie algebra.

The system of ODEs can be read--off directly  using Theorem \ref{prop_paul} and
setting 
$\Lambda=-\Theta(Y, Z, Z', -Y')$.
This
yields $\Lambda=(Y'Z-YZ')^{-1}$ and (\ref{pt1}) gives 
\be
\label{ode_tod}
Y''=\frac{2Y}{(Y'Z-YZ')^2}, \quad Z''=\frac{2Z}{(Y'Z-YZ')^2}.
\ee
The integral curves  are 
\[
\yy=A e^{\gamma X}+B e^{-\gamma X}, \quad \zz=
C e^{\gamma X}+D e^{-\gamma X},
\]
where $(A, B, C, D)$ are constants of integrations, and
$\gamma^2=(AD-BC)^{-1}\sqrt{2}^{-1}$. The original metric
can be recovered from the twistor lines by setting
\[
w=A+B, \quad z=(C+D), \quad y=\gamma(B-A), \quad x=\gamma (C-D).
\]
\subsection{Symmetry algebra of dimension 4}
Consider the system
\be
\label{ode_sym_4}
{Y}'' = 0, \quad Z'' = -\left( Z' + \sqrt{({Y'})^2 -1} \right)^2. 
\ee
This example was found by constructing the most general system of ODEs with 
Lie point symmetry algebra 
\begin{equation*} 
{L}_4 = \text{span} \{ \partial_X , \partial_Y , \partial_Z , Y \partial_X +
X \partial_Y \}
\end{equation*}
and imposing the torsion-free conditions. The algebra $L_4$ above is 
a particular realisation of the abstract algebra\footnote{The Referee
has pointed out that there are two more nonequivalent representations of 
this algebra: 
\[
L_{4a}= \text{span}\{\p_Y, -X\p_Y, \p_Z, \p_X-XZ\p_Y\}
\]
and
\[
L_{4b}= \text{span}\{2\p_X, \p_Z, -Y^2\p_X-y\p_Z, 2Z\p_X+\p_Y\}.
\]
We have checked that there is no torsion--free system of ODEs with symmetry
$L_{4b}$, and there exists a torsion--free system of ODEs with a seven--dimensional symmetry algebra which contains $L_{4a}$ as a sub-algebra.}
$A_{4,1}$ 
as given by Patera and Winternitz \cite{winter}. 
The integral curves of the  system (\ref{ode_sym_4}) are given by
\[
\yy=w+Xy, \quad \zz= \log{(X-x)}-\sqrt{y^2-1}X+z.
\]
Following the procedure applied in case of the system (\ref{boris})
we find the ASD conformal structure 
\[
g=dxdy+(dw+xdy)(dz+y\sqrt{y^2-1}^{-1} dw).
\]
This admits a null Killing vector $\p/\p z$ and thus fits into the classification of \cite{DW}.
\subsection{Symmetry algebra of dimension 3 or less}
If the symmetry algebra has dimension 3, then it necessarily belongs to the 
Bianchi classification of three--dimensional Lie algebras. There are many examples
in this case which are analytic continuations of Riemannian metrics.
See \cite{Tod} for a discussion of these examples. Examples of an ASD conformal
class with one or two symmetries can be found
in the Gibbons--Hawking class (\ref{GHmetric111}).
An axisymmetric  solution
$H$ to the wave equation (\ref{wave111}) gives a metric with two Killing vectors.
Another example with two--dimensional symmetry is the system (\ref{2Dsym})
discussed in the next Section.
A general solution of (\ref{wave111}) with no symmetries gives a metric  which 
only admits one Killing vector $\p/\p z$. Gravitational instantons of class
$D_k$ are examples of ASD conformal structures with no symmetries.

\section{Finsler structures with scalar flag curvature}
\label{section_finsler}
For a given $n$-dimensional manifold ${\z}$ with coordinates $X^k$ ($k = 1, \ldots, n$), a \emph{Finsler metric} is a positive continuous function $\mathcal{F}: \z \rightarrow [0,\infty)$ such that
\begin{itemize}
\item $\mathcal{F}$ is smooth on $T\z \backslash 0 = \{ (X^k,P^k) \in 
T\z | P \neq 0 \}$,
\item $\mathcal{F}(X^k,c P^k) = c \mathcal{F}(X^k,P^k)$ for $c > 0$,
\item The tensor $f_{ij} = \frac{1}{2} \frac{\partial^2 \mathcal{F}^2}{\partial P^i \partial P^j}$ is positive definite for all $(X^k,P^k) \in T\z \backslash 0$.
\end{itemize}
We shall consider the case $n=3$ and set  $X^k = (X,Y,Z)$ and $P^k = \dot{X}^k = \frac{dX^k}{dt}$ for some parameter $t$. Finsler geometry generalizes the notion of Riemannian geometry in that the norm on each tangent space $\mathcal{F}(X^k,\cdot)$ is not necessarily induced by a metric tensor. This makes these metrics useful in the study of problems involving paths of least time. The metric tensor $f_{ij}$ allows us to define Finslerian geodesics, which are integral curves of the system
\begin{equation*}
\ddot{X}^i + \gamma^i_{jk} \dot{X}^j \dot{X}^k = 0, \quad i,j,k = 1, \ldots, n
\end{equation*}
where
\begin{equation*}
\gamma^i_{jk} = \frac{1}{2} f^{il} \left( f_{lk, X^j} + f_{jl, X^k} - f_{jk, X^l} \right).
\end{equation*}
It was argued in \cite{crampsaund} that, for $n > 2$, given any system of ODEs with vanishing Wilczynski  invariants its integral curves arise as the set of unparametrised geodesics of a Finsler function of scalar flag curvature.
In this context, the torsion-free path geometries are viewed as projective equivalence classes of \emph{isotropic} sprays on $T\z$. Given a spray
\begin{equation*}
S = P^i \frac{\partial}{\partial X^i} - 2 \Gamma^i \frac{\partial}{\partial P^i} \,\,\,\,\,,\,\,\,\,\, i = 1, \ldots, n
\end{equation*}
we can define its Riemann curvature by
\begin{equation*}
R^l_{kij} = H_i \left( \Gamma^l_{jk} \right) - H_j \left( \Gamma^l_{ik} \right) + \Gamma^l_{im} \Gamma^m_{jk} - \Gamma^l_{jm} \Gamma^m_{ik}
\end{equation*}
where $\Gamma^i_{jk} = \frac{\partial^2 \Gamma^i}{\partial P^j \partial P^k}$, etc and the $H_i$ form a horizontal distribution determined by $S$:
\begin{equation*}
H_i = \frac{\partial}{\partial X^i} - \Gamma^j_i \frac{\partial}{\partial P^j}.
\end{equation*}
Then the spray $S$ is said to be \emph{isotropic} if its Jacobi endomorphism is given by
\begin{equation*}
R^i_j = R^i_{kjl} P^k P^l = \rho \delta^i_j + \tau_j P^i
\end{equation*}
for some function $\rho$ and covector $\tau_i$. It was shown in \cite{crampin}, that for $n > 2$, a spray is isotropic if and only if it is projectively equivalent to one whose Riemann curvature vanishes (so-called R-flat) and in \cite{cramptwo}, that any R-flat spray arises as the geodesic spray of some Finsler function. More importantly, in \cite{shen} it is shown that the geodesic sprays of a Finsler function are isotropic if and only if the Finsler function has scalar flag curvature. Here, the flag curvature of a Finsler metric is defined in terms of its Riemann curvature by
\begin{equation*}
K(X^k,P^k,V^k) = \frac{V^i (P^j R_{jikl} P^l) V^k}{f(P,P) f(V,V) - (f(P,V))^2}
\end{equation*}
where $V$ is a section of the tangent bundle transverse to $P$ and the indices are raised/lowered with the tensor $f_{ij}$. The flag curvature is scalar if
\begin{equation*}
K(X^k,P^k,V^k) = K(X^k,P^k).
\end{equation*}
Thus, we could construct all systems of ODEs with vanishing Wilczynski invariants if we know how to characterize the Finsler functions of scalar flag curvature. Although a lot of work has been done in this area, a complete characterization of such functions has not yet been achieved. Most success has come with a special type of Finsler function
\begin{equation*}
 \mathcal{F} = \sqrt{a_{ij} (X^k) P^i P^j} + b_i (X^k) P^i \,\,\,\,\,,\,\,\,\,\, i,j,k = 1, \ldots, n.
\end{equation*}
This is known as a \emph{Randers metric}. If we restrict our attention to Randers metrics of constant flag curvature, then we can use the following 
\begin{theo}[Bao--Robels--Shen \cite{bao}.]
A Randers metric $\mathcal{F}$ has constant flag curvature if and only if the corresponding Zermelo data $(h,W)$ satisfy the following:
\begin{itemize}
\item $h$ is a Riemannian metric with constant sectional curvature.
\item $W$ is a Killing vector or homothety of $h$.
\end{itemize}
\end{theo}
Here, the Randers data can be expressed in terms of the Zermelo data as follows \cite{robles}:
\begin{equation*}
a_{ij} = \frac{\lambda h_{ij} + W_i W_j}{\lambda^2} \,\,\,\,\,,\,\,\,\,\, b_i = - \frac{W_i}{\lambda}
\end{equation*}
where $\lambda = 1 - h_{ij} W^i W^j$ and $W_i = h_{ij} W^j$. This gives rise to a procedure for constructing torsion-free systems of ODEs from a three-dimensional Riemannian metric of constant sectional curvature and a homothety of this metric. The geodesic spray coefficients of such systems was worked out in \cite{robles}.
\\ \\ {\bf Example 4.}
Consider the Zermelo data
\begin{equation*}
h = dX^2 + dY^2 + Y^2 dZ^2, \qquad W = \frac{\partial}{\partial Z}.
\end{equation*}
The geodesics of the Randers metric associated to this Zermelo data are the integral curves of the systems of ODEs with vanishing Wilczynski invariants
\begin{eqnarray}
\label{2Dsym}
Y'' &=& \frac{2Y Z' \sqrt{1+ (Y')^2+ Y^2 \left((Z')^2-(Y')^2-1\right)}}
{Y (Y^2-1)}\\
&&+ \frac{Y\left(1+ (Y')^2+ (Z')^2 Y^2 ((Z')^2-(Y')^2-1)\right)}{(Y^2-1)^2},\nonumber \\
Z'' &=& \frac{2Y' \left( Z' + \sqrt{1+(Y')^2+ Y^2 \left((Z')^2 - (Y')^2 - 1\right)}\right)}{Y (Y^2-1)}.\nonumber
\end{eqnarray}
This system has two symmetries $\frac{\partial}{\partial X}$ and $\frac{\partial}{\partial Z}$.
\\\\ We can also view this correspondence in the other direction i.e, given a system of ODEs with vanishing Wilczynski invariants, we can construct a Finsler metric of scalar flag curvature.
\\ \\
{\bf Example 5.}
 To illustrate this point, let us consider the submaximal system
(\ref{submax}) corresponding to ASD Ricci--flat pp waves with constant
Weyl curvature.
Then, using the procedure in \cite{muzsnay}, this system describes the unparametrised geodesics of the Finsler function
\begin{equation*}
\mathcal{F} = \dot{X} {\mathcal G} \left( \frac{\dot{Y}}{\dot{X}} , 2\frac{\dot{Z}}{\dot{X}} -2 \frac{X \dot{Y}^3}{\dot{X}^3} + \frac{6Y \dot{Y}^2}{\dot{X}^2} \right) 
\end{equation*}
in some open domain of $\z$ where $f_{ij}$ is positive definite, where ${\mathcal G}$ is any function of two variables. In particular,  let ${\mathcal G}(x,y) = \sqrt{xy}$ 
and consider the  Lagrangian $L = \frac{1}{2} {\mathcal F}^2$ in unparametrised form
where $\dot{X}=1$
\begin{equation*}
L = 2Y' Z'-2X(Y')^4+6Y (Y')^3.
\end{equation*}
The Euler-Lagrange equations give (\ref{submax}).  For this example, the flag curvature of ${\mathcal F}$ vanishes.
\\ \\
Similarly other systems  with vanishing Wilczynski invariants arise from a variational principle induced by a Finsler structure.
Thus all these systems fit into the formalism of \cite{anderson}.
\section{Conclusions and outlook}
We have constructed several examples of systems of two second order ODEs
with vanishing Wilczynski invariants. The integral curves of such systems
are twistor curves which correspond to points in a four--manifold
$M$ with ASD conformal structure. The twistor curves can also be viewed as unparametrised
geodesics of Finsler structures with scalar flag curvature.
 
  The maximally symmetric, non--trivial example is
a Ricci--flat ASD $pp$--wave with 9--dimensional group of conformal symmetries,
and we have given examples of systems with symmetry groups
of dimensions between $9$ and $5$. Some of the examples
have a special form $Y''=0, \;\; Z''=G(X, Y, Z, Y', Z')$.
This is what Grossman calls 
`a weaker form of integrability for the second ruling' \cite{Gro}.
It corresponds to  the existence of a two parameter family
of $\beta$--surfaces (these are null ASD surfaces). This family gives 
a surface in the twistor space $\z$, and it is known \cite{DW}
any ASD $(2, 2)$ conformal structure with a conformal null Killing 
vector admits such structure. 
\vskip 5pt
In the case of trivial ODE with integral curves being straight 
lines in $\z=\R^3$
is a starting point for John's integral transform \cite{john}. Given
a function $\phi:\z\rightarrow \R$ with appropriate decay conditions at infinity
and a straight line $L\subset \z$, define a function $\hat{\phi}$
on the space $M$ of straight lines 
\[
\hat{\phi}(L)=\int_L\phi.
\]
The result of John is that the range
of this transform is characterised by an ultra-hyperbolic wave equation
on $M$, where the wave operator is induced by a flat $(2, 2)$ metric on $M=\R^4$.
If the line $L$ is parametrised by $X\rightarrow (X, Y=w+Xy, Z=z-Xx)$, then
\[
\hat{\phi}(w, z, x, y)=\int_\R \phi(X, w+Xy, z-Xx)dX\quad
\mbox{and} \quad \hat{\phi}_{wx}+\hat{\phi}_{zy}=0.
\]
It may be interesting to develop a non--linear version of John's transform
applicable to non--trivial systems of ODEs with vanishing Wilczynski invariants.
Integral curves of (\ref{submax}) could be a good starting point for this construction.


\end{document}